\documentclass[a4paper]{article}

\usepackage[dvips]{epsfig}
 \usepackage{graphics}
 \usepackage{graphicx}
\usepackage{color}
\usepackage[tight,footnotesize]{subfigure}
\usepackage{amsmath, amssymb, bbm,vmargin}

\begin{document}



\title{Predictive maintenance for the heated hold-up tank}


\author{Beno\^\i te de Saporta \and Huilong Zhang}
\date{}
\maketitle

\begin{abstract}
We present a numerical method to compute an optimal maintenance date for the test case of the heated hold-up tank. The system consists of a tank containing a fluid whose level is controlled by three components: two inlet pumps and one outlet valve. A thermal power source heats up the fluid. The failure rates of the components depends on the temperature, the position of the three components monitors the liquid level in the tank and the liquid level determines the temperature. Therefore, this system can be modeled by a hybrid process where the discrete (components) and continuous (level, temperature) parts interact in a closed loop. 
We model the system by a piecewise deterministic Markov process, propose and implement a numerical method to compute the optimal maintenance date to repair the components before the total failure of the system. 
\end{abstract}

\section{Introduction}
\label{section:intro}

A complex system is inherently sensitive to failures of its components. One must therefore determine maintenance policies in order to maintain an acceptable operating condition. Optimizing the maintenance is a very important problem in the analysis of complex systems. It determines when it is best that maintenance tasks should be performed on the system in order to optimize a cost function: either maximize a performance function or conversely minimize a loss function. Moreover, this optimization must take into account the random nature of failures and random evolution and dynamics of the system. 

The example considered here is the maintenance of the heated hold-up tank, a well know test case for dynamic reliability, see e.g. \cite{devooght97, MZ95, MZ96, ZDDG08}. The system consists of a tank containing a fluid whose level is controlled by three components: two inlet pumps and one outlet valve. A thermal power source heats up the fluid. The failure rate of the components depends on the temperature, the position of the three components monitors the liquid level in the tank, and in turn, the liquid level determines the temperature. The main characteristic of this system is that it can be modeled by a stochastic hybrid process, where the discrete and continuous parts interact in a closed loop. As a consequence, simulating this process and computing related reliability indices has been a challenge for the dynamic reliability community. To our best knowledge, optimization of maintenance policies for the heated hold-up tank has not been addressed yet in the literature.  

The only maintenance operation considered here is the complete replacement of all the failed components and the system restarts in its initial equilibrium state. Partial repairs are not allowed. Mathematically, this problem of preventive maintenance corresponds to a stochastic optimal stopping problem as explained by example in the book of Aven and Jensen \cite{AJ99}. It is a difficult problem because of the closed loop interactions between the state of the components and the liquid level and temperature. A classical approach consists in using condition-based maintenance (CBM) to act on the system based on its current state and before its failure. One can for example calculate the remaining useful life (RUL) of the system and the preventive replacement is carried out when the deterioration level exceeds a certain threshold or enters in a certain state \cite{vannoortwijk, grall}. Our approach also takes into account the current state of the process, but our decision rule is not based on damage accumulation nor does it correspond to hitting some threshold. Instead, it involves a performance function that reflects that the longer the system is in a functioning state the better.

The dynamics of the heated hold-up tank can be modeled by a piecewise deterministic Markov process  (PDMP), see \cite{ZDDG08}. Therefore, our maintenance problem boils down to an optimal stopping problem for PDMP's. PDMP's are a class of stochastic hybrid processes that has been introduced by Davis \cite{Davis93} in the 80's. These processes have two components: a Euclidean component that represents the physical system (e.g. temperature, pressure, \ldots) and a discrete component that describes its regime of operation and/or its environment. Starting from a state $x$ and mode $m$ at the initial time, the process follows a deterministic trajectory given by the laws of physics until a jump time that can be either random (e.g. it corresponds to a component failure or a change of environment) or deterministic (when a magnitude reaches a certain physical threshold, for example the pressure reaches a critical value that triggers a valve). The process restarts from a new state and a new mode of operation, and so on. This defines a Markov process. Such processes can naturally take into account the dynamic and uncertain aspects of the evolution of the system. A subclass of these processes has been introduced by Devooght \cite{devooght97} for an application in the nuclear field. The general model has been introduced in dynamic reliability by Dutuit and Dufour \cite{DD02}.


The objective and originality of this paper is twofold. First, we propose an optimization procedure for a well-known benchmark in the dynamic reliability literature. The tank model was first introduced by \cite{aldemir87a} where only one continuous variable (liquid level)
 is taken into account, and then in \cite{marseguerra1994a} and
\cite{MZ95} where the second variable (temperature) is introduced. They have tested various Monte
Carlo approaches to simulate the process to compute reliability and safety indices. In \cite{tombuyses1996a}, the authors have used the
same system to present continuous cell-to-cell mapping Markovian
approach (CCCMT) still to simulate the process. The simulation of the holdup tank example has been and is still widely studied in the literature (not exhaustive)
\cite{siu1994a,cojazzi1996a,dutuit1997a,schoenig2006a,li2011a,JRR12}. Here we go one step further and not only propose to simulate the tank process but also we optimize it.

Second, even though PDMP's have been recognized as a powerful modeling tool for dynamic reliability problems \cite{devooght97,DD02}, there are very few numerical tools adapted to these processes. Our aim is to further demonstrate the high practical power of the theoretical methodology described in \cite{AAP10}, by applying it to the tank benchmark. In \cite{AAP10}, the authors have proposed a numerical algorithm to optimize PDMP's and have studied its theoretical properties.
This optimization procedure was first applied to an example of maintenance of a metallic structure subject to corrosion, without closed loop interactions or deterministic jumps. In addition, the system has only one continuous variable and the cost function is simple and does not depend on time, see \cite{JRR12}.
In this paper, we adapt the numerical procedure proposed in \cite{AAP10} to the more challenging heated hold-up tank problem with two continuous variables, deterministic jumps when these variables hit some given boundaries and closed loop interactions between continuous and discrete variables. Furthermore, we consider a cost function that depends on both continuous variables as well as on the running time.

The remainder of this paper is organized as follows. In section \ref{section:model}, the dynamics of the heated hold-up tank is presented with more details as well as the framework of PDMP's. In section \ref{section:opti} the formulation of the optimal stopping problem for PDMP's and its theoretical solution are briefly recalled and the four main steps of the algorithm are detailed. In section \ref{section:res} the numerical results obtained on the example of the tank are presented and discussed. Finally, in section \ref{section:conclusion} a conclusion and perspectives are presented.

\section{Model}
\label{section:model}
We are interested in the maintenance of a heated hold-up tank. The dynamics of the tank can be modeled by a piecewise deterministic Markov process  (PDMP). We first describe with more details the dynamics of the tank, then we recall the definition and some basic properties of PDMP's. The tank model is a well known benchmark in dynamic reliability. It was first introduced by \cite{aldemir87a} where only one continuous variable (liquid level) is taken into account, and then in \cite{marseguerra1994a} and \cite{MZ95} where the second variable (temperature) is introduced. We have kept the values of the parameters defined in those papers.
\subsection{The heated hold-up tank}
\label{section:tank}
\begin{figure}[htb]
\begin{center}
\includegraphics[width=8cm]{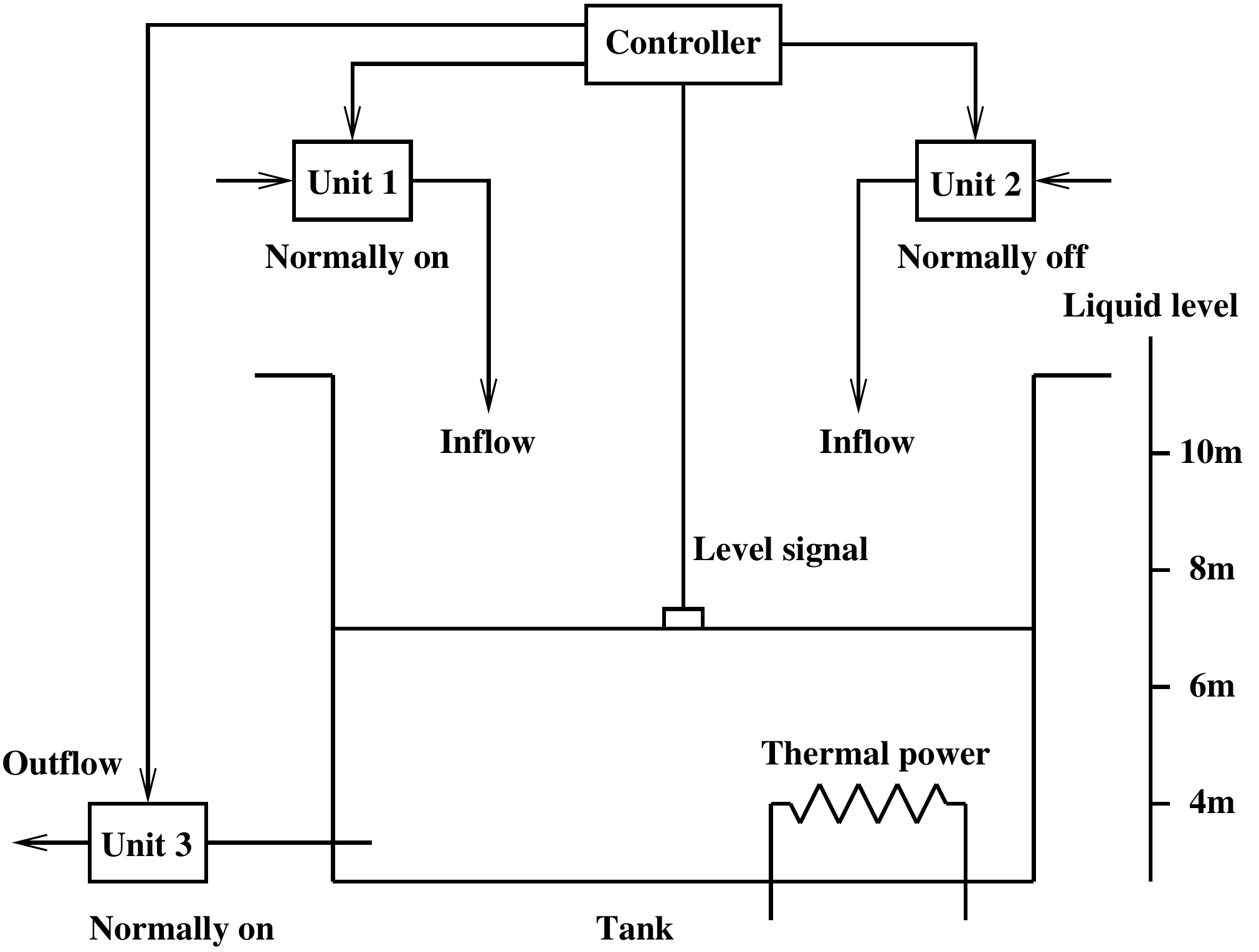}
\caption{The heated hold-up tank}\label{fig1}
\end{center}
\end{figure}

\begin{figure}[htb]
\begin{center}
\includegraphics[width=6cm]{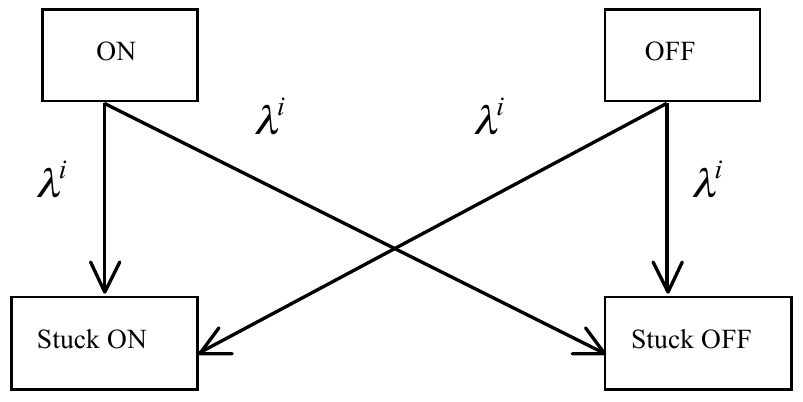}
\caption{Transitions for unit $i$}\label{fig2}
\end{center}
\end{figure}

The system is represented on Figure \ref{fig1}. It consists of a tank containing a fluid whose level is controlled by three components: two inlet pumps (units 1 and 2) and one outlet valve (unit 3). A thermal power source heats up the fluid. The variables of interest are the liquid level $h$, the liquid temperature $\theta$ and the state of the three components and the controller. 
Each component has four states: ON, OFF, Stuck ON or Stuck OFF. Once a unit is stuck (either on or off) it cannot change state. The possible transitions between these four states are given in Figure \ref{fig2}. Thus, by a random transition a working unit can only become stuck (either on or off). The initial state of the components is ON for units 1 and 3 and OFF for unit 2.
The intensity of jumps $\lambda^i$ for unit $i$ depends on the temperature through the equation $\lambda^i=a(\theta)l^i$ with $a(\theta)$ given in Eq.~(\ref{eq:a}), see \cite{marseguerra1994a,MZ95}
\begin{equation}\label{eq:a}
a(\theta)=\frac{b_1\exp\big(b_c(\theta-20)\big)+b_2\exp\big(-b_d(\theta-20)\big)}{b_1+b_2}.
\end{equation}
 
\begin{figure}[t]
\begin{center}
\includegraphics[width=5cm]{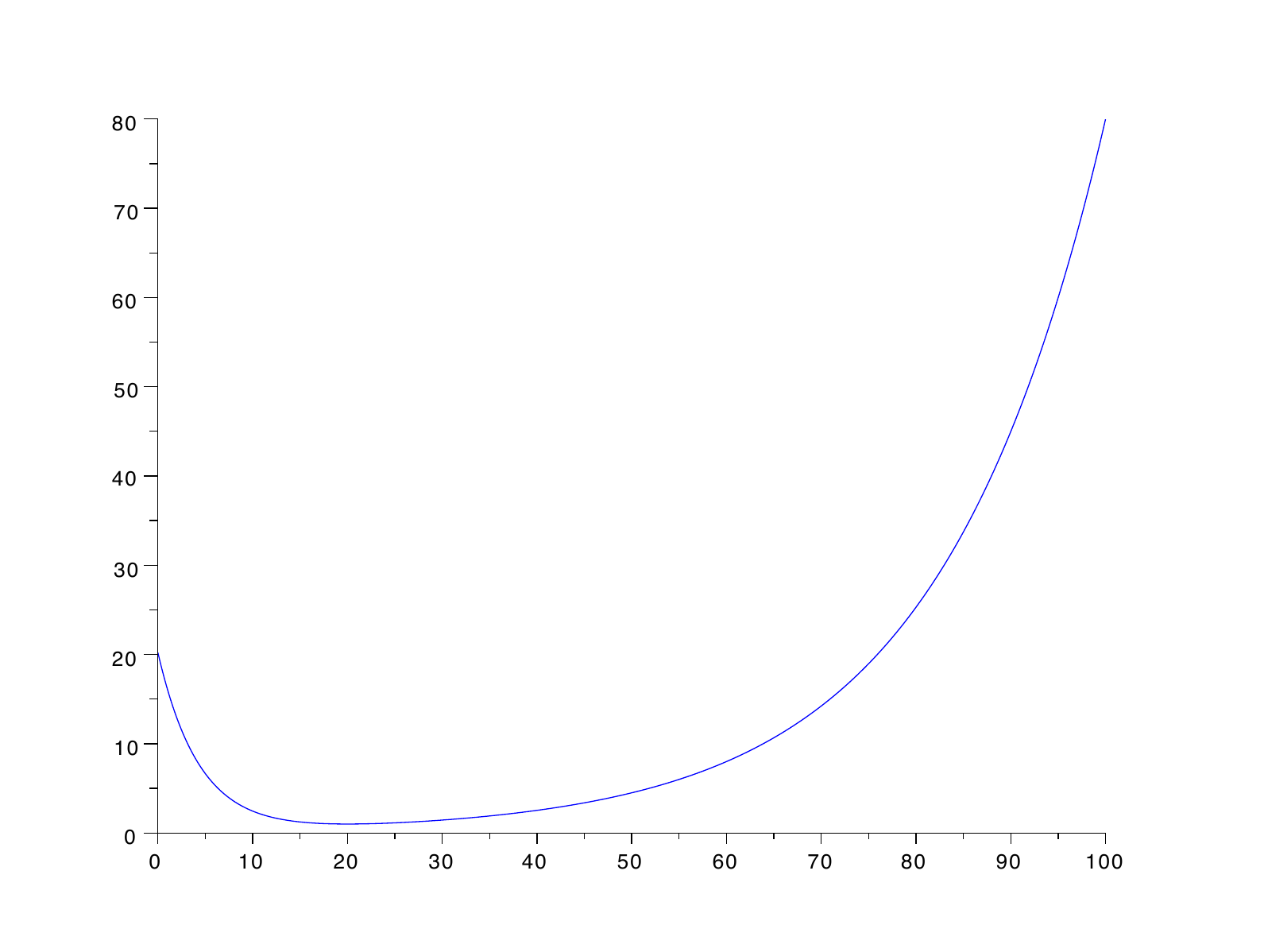}
\caption{$a(\theta)$ as a function of $\theta$}\label{fig3}
\end{center}
\end{figure}

Function $a(\theta)$ is represented on Figure \ref{fig3} and the various parameters come from the literature, see \cite{marseguerra1994a,MZ95}, and are given in Table \ref{tab1}. 
The special form of the failure rate $\lambda^i$ as a product of a constant depending on $i$ and a function of the temperature allows for all three units to have failure rates with the same dependence on the temperature, but different scaling parameters. Indeed, at the reference temperature of $20^\circ C$, the mean time to failure of unit $1$ is $438h$, for unit $2$, it is $350h$ and for unit $3$ it is $640h$.

In addition,
the shape for function $a(\theta)$ was chosen in the original benchmark so that there is a very high failure rate when the temperature is high. More specifically, the parameters are chosen such that $a(\theta)$ is lowest (equal to $1$) when the temperature is equal to a reference temperature of $20^\circ C$, it equals $20$ when the temperature is $0^\circ C$ and it is highest (equal to $80$) when the temperature equals the critical temperature of $100^\circ C$. The exponential functions are chosen in order to enable this very high dependence with the temperature. Roughly speaking, the units fail $80$ times more often when the temperate is $100^\circ C$ than when it is $20^\circ C$.

In addition, control laws are used to modify the state of the components to keep the liquid within two acceptable limits: $6$ meters and $8$ meters. If the liquid level drops under $6\ m$, the components 1, 2, 3 are put respectively in state ON, ON and OFF (provided they are not stuck). If the liquid level rises above $8\ m$, the components are put respectively in the state OFF, OFF and ON (provided they are not stuck). Unlike the classical model presented in \cite{devooght97, MZ95, MZ96, ZDDG08}, we also allow the control unit to fail. At each solicitation, the control may succeed with probability $p=0.8$ independently from previous successes. Once it has failed, it will never succeed again. Therefore the control unit has two possible states: working $1$ or failed $0$.
\begin{table}
\begin{center}
\begin{tabular}{|c|c|}
\hline
Parameter&Value\\
\hline
$b_1$&$3.0295$\\
$b_2$&$0.7578$\\
$b_c$&$0.05756$\\
$b_d$&$0.2301$\\
$\theta_{in}$&$15^\circ C$\\
$l^1$&$2.2831\cdot10^{-3}\  h^{-1}$\\
$l^2$&$2.8571\cdot10^{-3}\  h^{-1}$\\
$l^3$&$1.5625\cdot10^{-3}\  h^{-1}$\\
$G$&$1.5\  mh^{-1}$\\
$K$&$23.88915\ m^\circ Ch^{-1}$\\
\hline
\end{tabular}
\caption{Parameters for the tank dynamics}\label{tab1}
\end{center}
\end{table}

The evolution of the liquid level $h$ depends on the position of the three components through the differential equation~(\ref{eq:h})
\begin{equation}\label{eq:h}
\frac{\partial h}{\partial t}=(\nu_1+\nu_2-\nu_3)G,
\end{equation}
where $\nu_i=1$ if component $i$ is ON or Stuck ON, $\nu_i=0$ otherwise, and $G$ is the common flow of the three components and is given in Table \ref{tab1}. The initial level is $h_0=7\ m$. Eq.~(\ref{eq:h}) simply means that each pump on contributes to rise the liquid level, whereas if the outlet valve is on, it contributes to the decrease of the liquid level.
The temperature $\theta$ depends on the liquid level through the differential equation~(\ref{eq:theta})
\begin{equation}\label{eq:theta}
\frac{\partial \theta}{\partial t}=\big((\nu_1+\nu_2)G(\theta_{in}-\theta)+K\big)h^{-1},
\end{equation}
where $\theta_{in}$ is the temperature of the incoming fluid, and $K$ is a constant of the tank, the values of these parameters are given in Table \ref{tab1}. As the tank is heated and the incoming liquid has a constant temperature $\theta_{in}$, Eq.~(\ref{eq:theta}) reflects that the temperature converges to an equilibrium state as long as there is some incoming fluid. The temperature can increase to the threshold $100^\circ C$ if there is no incoming fluid. The initial temperature is $\theta_0=30.9261^\circ C$, so that the system is initially at an equilibrium state, and nothing happens until the failure of one of the components. The system stops as soon as one of the top events is reached: dry out ($h<4\ m$), overflow ($h>10\ m$) or hot temperature ($\theta>100^\circ C$).
\subsection{Piecewise deterministic Markov processes}
\label{section:PDP}
As in \cite{ZDDG08}, we model the tank by a Piecewise-deterministic Markov processes (PDMP). PDMP's are a general class of hybrid processes. They are defined as follows. Let $M$ be the finite set of the possible modes of the system. In our tank example, the modes correspond to the possible positions of the inlet pumps, outlet valve and control unit. The components can be ON, OFF, stuck ON or Stuck OFF, the control unit can be in position $0$ or $1$. Therefore, there are 128 possible modes for our system (but only 74 can actually be reached from the equilibrium starting point). 

For all modes $m$ in $M$, let $E_m$ be an open subset in $\mathbb{R}^d$.  In our case $d=3$ as we need to take into account the running time as well as the liquid level and temperature, and $E_m$ is a subset of $(4,10)\times[15,100)\times[0,+\infty)$ (depending on $m$).
A PDMP is 
defined from three local characteristics $(\Phi, \lambda, Q)$ where
\begin{itemize} 
\item the flow $\Phi : M\times \mathbb{R}^d\times \mathbb{R} \rightarrow \mathbb{R}^d$ is continuous 
and for all $ s, t \geq 0$,  one has $\Phi(\cdot, \cdot,t+s) = \Phi(\Phi(\cdot, \cdot,s),t)$. 
It describes 
the deterministic trajectory of the process between jumps. In the tank example, it is given by the solution of Eq.~(\ref{eq:h}) and (\ref{eq:theta}).
For all $(m, \textbf{\em{x}})$ 
in $ M \times E_m$, set
\begin{equation} 
t^\ast(m,\textbf{\em{x}}) = \inf \{t>0 : \Phi(m,\textbf{\em{x}},t) \in \partial E_m\},
\end{equation} 
the time to reach the  boundary of the domain starting from $\textbf{\em{x}}$ in mode $m$. For the tank, the boundary is one of the thresholds $4\ m$, $6\ m$, $8\ m$, $10\ m$, $100^\circ C$ or $1000$ hours for the running time.
\item the jump intensity $\lambda$ characterizes the frequency of jumps. For 
all $(m, \textbf{\em{x}})$ in $M \times E_m$, and $t \leq t^\ast (m, \textbf{\em{x}})$, set
\begin{equation}
\Lambda(m,\textbf{\em{x}},t) = \int_0^t \lambda (\Phi(m,\textbf{\em{x}},s))\,ds.
\end{equation}  
For the tank the jump intensity given a mode $m$ is the sum of the intensities $\lambda_i$ for the remaining possible jumps of the three units.
\item the Markov kernel $Q$ represents the transition measure of the process 
and allows to select the new location and mode after each jump. In our example, $Q$ acts only on the mode components and leaves the liquid level $h$, temperature $\theta$ and running time unchanged. It selects one of the remaining possible failures of the three components, or corresponds to an attempted control law.
\end{itemize} 

\begin{figure}[htb]
\begin{center}
\includegraphics[width=12cm]{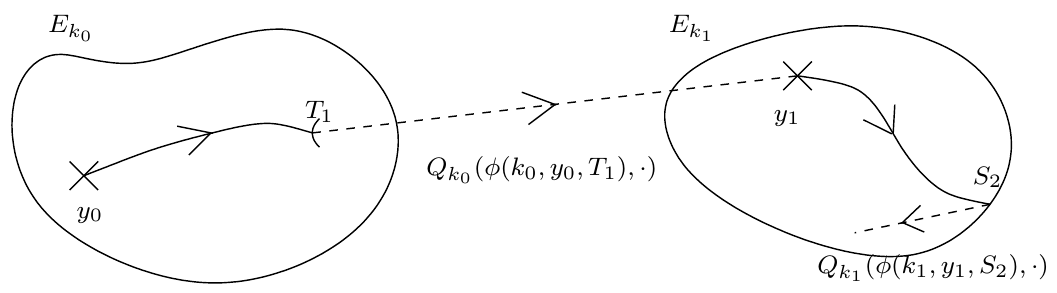}
\caption{An example of path for a PDMP until the second jump. The first jump is 
random. The second jump is deterministic because the process has reached the 
boundary of the domain.}\label{fig4}
\end{center}
\end{figure}

The trajectory $\textbf{\em{X}}_t = (m_t, \textbf{\em{x}}_t)$ of the process can then be defined 
iteratively. It starts at an initial point $\textbf{\em{X}}_0 = (k_0, y_0)$ with $ k_0\in M$
 and $y_0\in E_{k_0}$. For the tank, $k_0=(ON, OFF, ON, 1)$ and $y_0=(7, 30.9261,0)$. The first jump time $T_1$ is determined by Eq.~(\ref{eq:T})
\begin{equation}\label{eq:T}
\mathbb{P}_{(k_0,y_0)}(T_1>t) = \left \{ 
\begin{array} {lcl}
e ^{-\Lambda(k_0,y_0,t)} &\textrm{if} &~t<t^\ast(k_0,y_0),\\
0 &\textrm{if}  &~t\geq t^\ast(k_0,y_0).
\end{array} 
\right .
\end{equation} 
It corresponds to the first failure time of one of the components as in our case $t^\ast(k_0,y_0)=+\infty$.
On the interval $[0, T_1)$, the process follows the  deterministic trajectory 
$m_t = k_0$ and $\textbf{\em{x}}_t = \Phi(k_0, y_0, t)$. At the random time $T_1$, 
a jump occurs. Note that in general a jump can be either a discontinuity in the Euclidean variable $\textbf{\em{x}}_t$ or a change of mode. The process restarts at a new mode and/or position $\textbf{\em{X}}_{T_1}=(k_1,y_1)$, according 
to distribution $Q_{k_0}(\Phi(k_0, y_0, T_1), \cdot)$. An
inter jump time  $T_2- T_1$ is then selected in a similar way, and on the interval $[T_1, T_2)$ the process 
follows the path $m_t=k_1$ and $\textbf{\em{x}}_t = \Phi(k_1, y_1, t - T_1)$. 
Thereby, iteratively, a PDMP is constructed, see Figure \ref{fig4} for an 
illustration. 

Let $\textbf{\em{Z}}_0 = \textbf{\em{X}}_0$, and for $n \geq 1, ~\textbf{\em{Z}}_n = \textbf{\em{X}}_{T_N}$, location and mode of the 
process after each jump. Let $S_0 = 0$, $S_1=T_1$ and for $n \geq 2$, 
$S_n = T_n-T_{n-1}$ the inter-jump times between two consecutive jumps, 
then $(\textbf{\em{Z}}_n, S_n)$  is a Markov chain, which is the only source of 
randomness of the PDMP and contains all information on its random part.
 Indeed, if one knows the jump times and the positions after each jump, one can 
reconstruct the deterministic part of the trajectory between jumps. Namely, if one knows the time and nature of all the components failures, one can reconstruct the history of the liquid level and temperature though Eq.~(\ref{eq:h}) and (\ref{eq:theta}). It is a 
very important property of PDMP's that is at the basis of the numerical procedure.
\section{Optimization problem}
\label{section:opti}
\subsection{General framework}
\label{section:framework}
The general mathematical problem of optimal stopping corresponding to this maintenance problem can be found in \cite{gugerli86, AAP10, JRR12}. It is now briefly recalled. Let 
$\textbf{\em{z}} = (k_0, y_0)$ be the starting point of the PDMP $(\textbf{\em{X}}_t)$. Let $\mathcal{M}_f$ be the 
set of all stopping times $\tau$ for the natural filtration of the PDMP ($\textbf{\em{X}}_t$) satisfying 
$\tau \leq T_f$ that is to say that the intervention takes place before the time $T_f=1000\ h$. It has been shown in previous studies that by $1000\ h$, all the units are stuck and the events of interest have been observed, see e.g. \cite{devooght97, MZ95, MZ96, ZDDG08}.
 Let $g$
 be the cost function to optimize. Here, $g$ is a reward function that has to be maximized. The optimization problem to solve is given in Eq.~(\ref{eq:v})
\begin{equation}\label{eq:v}
v(\textbf{\em{z}}) = \sup_{\tau \in M_f}E_{\textbf{\em{z}}}\left [g(\textbf{\em{X}}_\tau)\right ].
\end{equation}
The function $v$ is called the {\em value function} of the problem and represents 
the maximum performance that can be achieved. Solving the optimal stopping 
problem is firstly to calculate the value function, and secondly to find a 
stopping time $\tau$ that achieves this maximum. This stopping time is 
important from the application point of view since it corresponds to the 
optimum time for maintenance. 
In general, such an optimal stopping time does not exist.  Define then 
$\epsilon$-optimal stopping times as achieving optimal value minus $\epsilon$, 
i.e. $v(\textbf{\em{z}})-\epsilon$.

Under fairly weak regularity conditions, Gugerli has shown in \cite{gugerli86} 
that the value function $v$ can be calculated iteratively as follows. 
First, choose the computational horizon $N$ such that after $N$ jumps, the running time $t$ has reached $T_f$  for almost all trajectories.
Let $v_N=g$ be the reward function, and iterate an operator $L$ backwards, see Eq.~(\ref{eq:systeme}). 
The function $v_0$ thus obtained is equal to the value function $v$.
\begin{equation} \label{eq:systeme}
\left \{
\begin{array} {lcl}
v_N & = & g, \\
v_k & = & L(v_{k+1},g), \quad0\leq k\leq N-1.
\end{array} 
\right .
\end{equation} 
Operator $L$ defined in Eq.~(\ref{def L}) is complex and involves a 
continuous maximization, conditional expectations and indicator functions, 
even if the reward function $g$ is very regular:
\begin{eqnarray} \label{def L}
\lefteqn{L(w,g)(\textbf{\em{z}})}\nonumber\\ 
&\equiv &  \sup_{u\leq t^\ast(\textbf{\em{z}}) }
\left \{ 
E \left [w(\textbf{\em{Z}}_1)\mathbbm{1}_{\{S_1<u\wedge t^\ast(\textbf{\em{z}})\}} 
+ g(\Phi(\textbf{\em{z}},u))\mathbbm{1}_{\{S_1\geq u \wedge t^\ast(\textbf{\em{z}})\}}|\textbf{\em{Z}}_0=\textbf{\em{z}} \right ] 
\right \}\nonumber\\
 && \vee E \left [ w(\textbf{\em{Z}}_1)|\textbf{\em{Z}}_0=\textbf{\em{z}} \right ].
\end{eqnarray} 
However, this operator depends only on the discrete time 
Markov chain $(\textbf{\em{Z}}_n, S_n)$. Gugerli also proposes an iterative construction of $\epsilon$-optimal 
stopping times, which is too technical to be described here, 
see \cite{gugerli86} for details.

\begin{figure}[htb]
\begin{center}
\includegraphics[width=10cm]{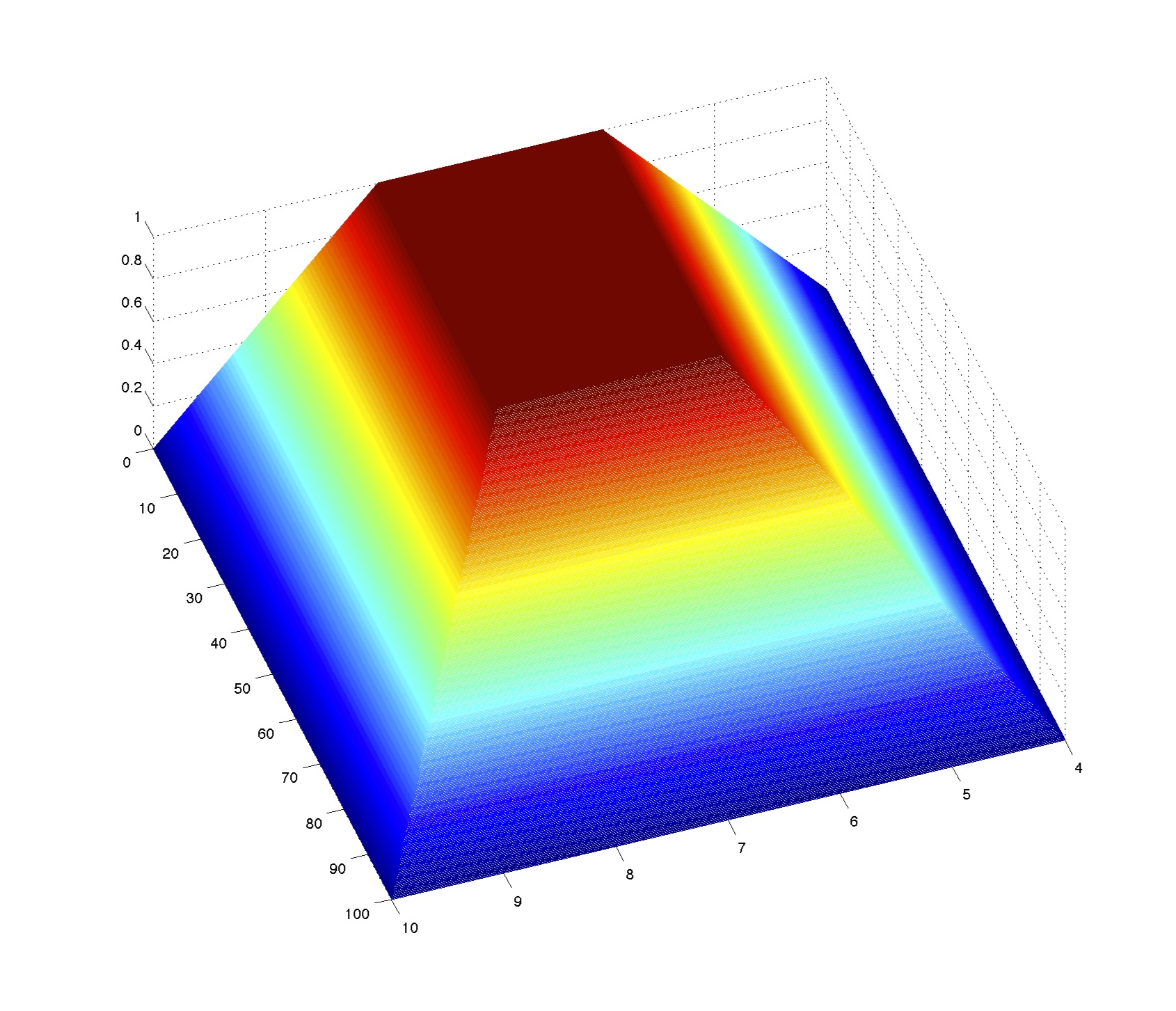}
\caption{Reward function $f$ as a function of $h$ and $\theta$}\label{fig5}
\end{center}
\end{figure}

In our example, the reward function has two components $g(h, \theta, t)=f(h, \theta)t^{\alpha}$. The first one $f$ depends on the liquid level and temperature, and reflects that the reward is maximal (set to $1$) when $h$ and $\theta$ are in the normal range ($6\ m \leq h\leq8\ m$, $\theta\leq 50^\circ C$), minimal (set to $0$) when reaching the top events: dry out ($h<4\ m$), overflow ($h>10\ m$) or hot temperature ($\theta >100^\circ C$) and continuous in between, see Figure \ref{fig5} for an illustration. The second term tα involves the time and reflects that the longer the system is functioning the higher the reward. The parameter $\alpha$ is set to 1.01 for smoothness.
%
\subsection{Numerical procedure}
\label{section:num proc}
In \cite{AAP10, JRR12} the authors propose a numerical method to approximate the value function for the optimal stopping problem of a general PDMP. The approach is based on a discretization of the operator L defined above. Our algorithm for calculating the value function is divided into three stages: a discretization of the Markov chain $(\textbf{\em{Z}}_n, S_n)$, a path-adapted time discretization between jumps, and finally a recursive computation of the value function~$v$. Then, the calculation of a quasi-optimal stopping time only uses comparisons of quantities already calculated in the approximation of the value function, which makes this technique particularly attractive. These stages are briefly recalled below.

\subsubsection{Quantization}
\begin{figure}[h]
\centering
\subfigure[Standard normal density in 2D]
{\includegraphics[width=0.49\linewidth]{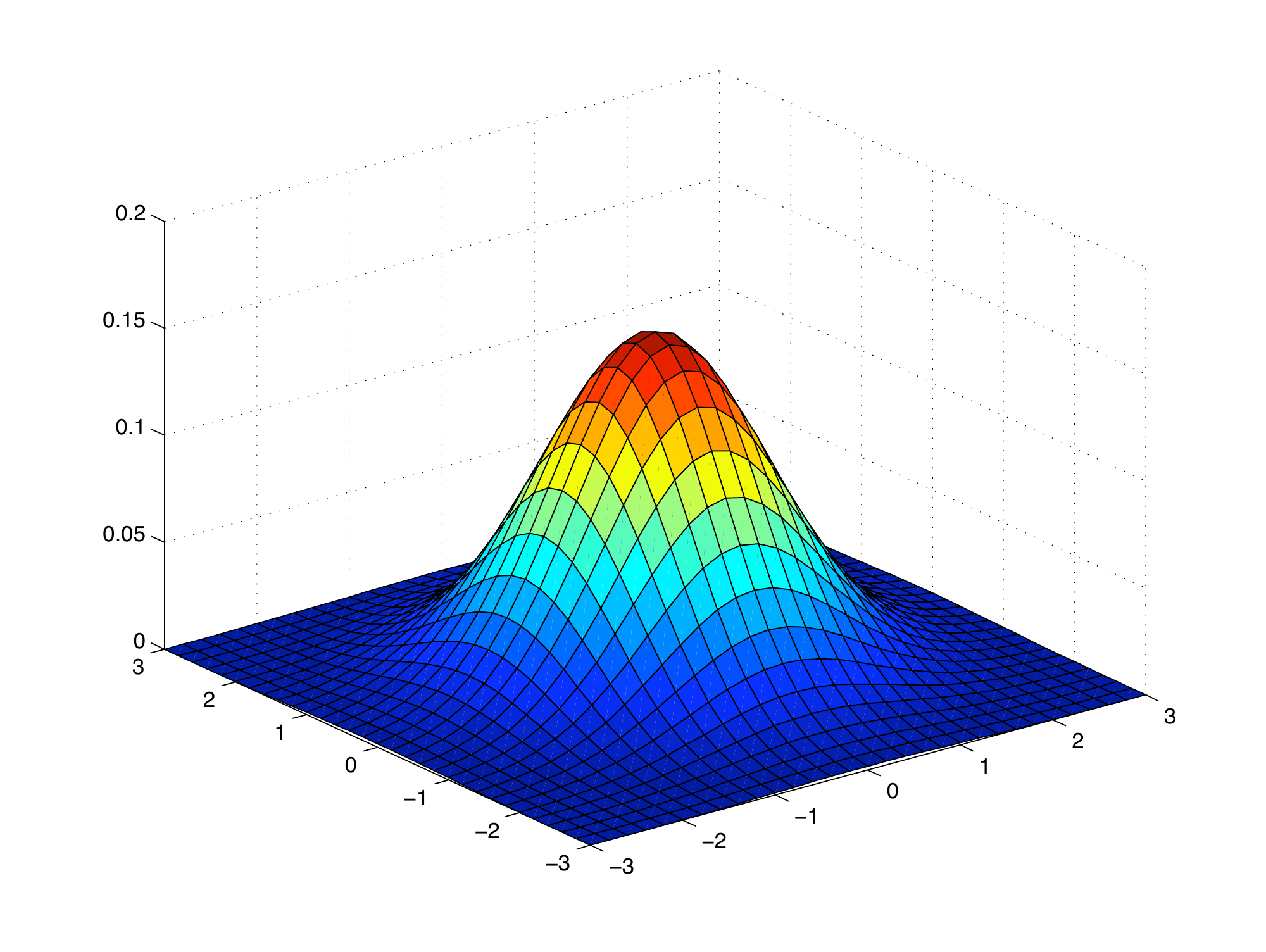}} 
\subfigure[Quantization grid]
{\includegraphics[width=0.49\linewidth]{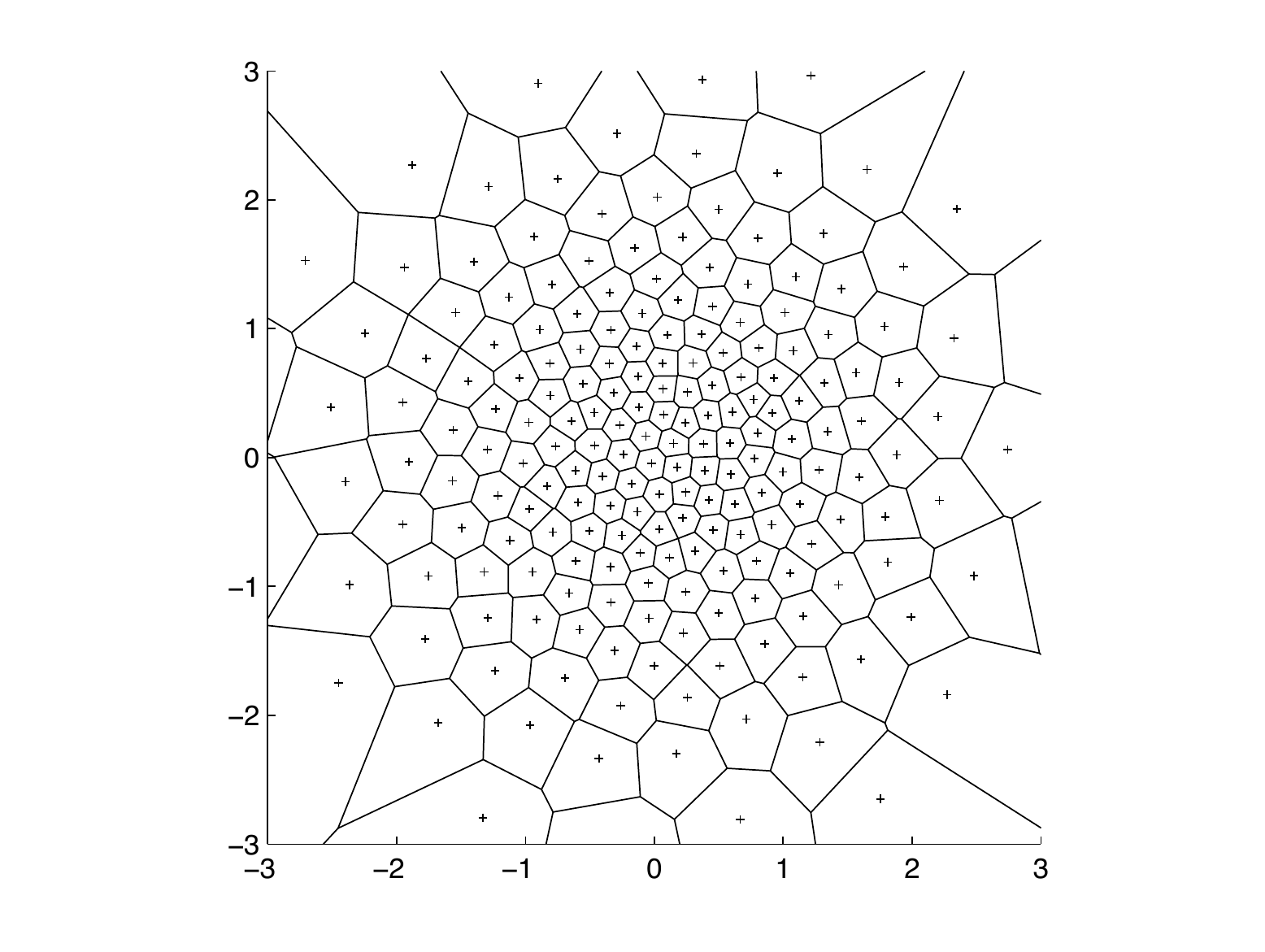}} 
\caption{Example of quantization grid for a normal distribution (200 points)}
\label{fig6}
\end{figure}

The goal of the first step is to approximate the continuous state space 
Markov chain $(\textbf{\em{Z}}_n, S_n)$ 
by a  discrete state space sequence 
$(\widehat{\textbf{\em{Z}}}_n, \widehat{S}_n)$. To this aim, we use the quantization algorithm described in 
details in e.g. \cite{pages98, pages05, pages04a, pages04b}.
Roughly speaking, more points are put in the areas of high density of the random variable, see Figure \ref{fig6} for an example of quantization grid for the standard
normal distribution in two dimensions. 
The quantization algorithm is based 
on Monte Carlo simulations combined with a stochastic gradient method. 
It provides $N+1$ grids, one for each couple $(\textbf{\em{Z}}_n,S_n)$ ($0\leq n\leq N)$, with a fixed number of points in each grid. The algorithm also
provides weights for the grid points and probability transition between two points of two
consecutive grids fully determining the distribution of the approximating sequence $(\widehat{\textbf{\em{Z}}}_n, \widehat{S}_n)$. 
The quantization theory ensures that the $L^2$ 
distance between $(\widehat{\textbf{\em{Z}}}_n, \widehat{S}_n)$ and $(\textbf{\em{Z}}_n, S_n)$ tends to 0 as the 
number of points in the  quantization grids tends to infinity, see  \cite{pages04a}.
%
\subsubsection{Time discretization}
Now the continuous maximization of the operator $L$ is replaced by a 
finite maximization, that is to say that one must discretize the time 
intervals $[0, t^\ast(\textbf{\em{z}})]$ for a finite number of $\textbf{\em{z}}$, namely each $\textbf{\em{z}}$ in the quantization grids. For this, choose time steps $\Delta(\textbf{\em{z}}) <t^\ast(\textbf{\em{z}})$ and take a regular discretization $G(\textbf{\em{z}})$ of $[0, t^\ast(\textbf{\em{z}})-\Delta(\textbf{\em{z}})]$ with step $\Delta(\textbf{\em{z}})$.
The maximum in such grids is less than $t^\ast(\textbf{\em{z}}) - \Delta(\textbf{\em{z}})$, which is a crucial 
property to derive error bounds for the algorithm, see~\cite{AAP10}. 
%
\subsubsection{Approximate calculation of the value function}
One now has all the tools to provide an approximation of the operator $L$ given in Eq.~(\ref{Lchap}). For 
each $1\leq n \leq N$, and for all $\textbf{\em{z}}$ in the quantization grid at time $n-1$, 
 set
\begin{eqnarray} \label{Lchap}
\lefteqn{ \widehat{L}_n (w,g)(\textbf{\em{z}}) }\nonumber\\
&\equiv& \max_{u\in G(\textbf{\em{z}}) }
\Big \{ 
E  \Big [ w(\widehat{\textbf{\em{Z}}}_{n-1})\mathbbm{1}_{\{\widehat{S}_{n}<u\wedge t^\ast(\textbf{\em{z}})\}} + g(\Phi(\widehat{\textbf{\em{Z}}}_{n-1},u))\mathbbm{1}_{\{\widehat{S}_n\geq u \wedge t^\ast(\textbf{\em{z}})\}} 
  |\widehat{\textbf{\em{Z}}}_{n-1}=\textbf{\em{z}} \Big ] \Big \}\nonumber\\
&& \qquad \vee E \left [ w(\widehat{\textbf{\em{Z}}}_n)|\widehat{\textbf{\em{Z}}}_{n-1}=\textbf{\em{z}} \right ].
\end{eqnarray} 
Note that because there are different quantized approximations at each time step, there also are different discretizations of operator $L$ at each time step. It should be also noted that the 
conditional expectations taken with respect to a process with finite state space 
are actually finite weighted sums.
One then constructs an approximation of the value function by backward iterations of the operators $ \widehat{L}_n$:
\begin{equation} 
\left \{
\begin{array} {lcl}
\widehat{v}_N & = & g, \\
\widehat{v}_{n-1} (\widehat{\textbf{\em{Z}}}_{n-1})& = & \widehat{L}_n(\widehat{v}_{n},g)(\widehat{\textbf{\em{Z}}}_{n-1}), 
\quad 1\leq n\leq N.
\end{array} 
\right .
\end{equation} 
Then take $\widehat{v}_0(\widehat{\textbf{\em{Z}}}_0)=\widehat{v}_0(\textbf{\em{z}})$ as an approximation of the 
value function $v$ at the starting point $\textbf{\em{z}}$ of the PDMP. The difference between $\widehat{v}_0(\textbf{\em{z}})$ and $v$ goes to zero as the number of points in the quantization grids goes to infinity, see
\cite{AAP10} for details and a convergence rate.
%
\subsubsection{Calculation of a quasi-optimal stopping time}
A method to compute an $\epsilon$-optimal stopping 
time has also been implemented. The discretization is much more complicated and subtle 
than that of operator $L$, because one needs both to use the true Markov chain 
$(\textbf{\em{Z}}_n, S_n)$ and its quantized version $(\widehat{\textbf{\em{Z}}}_n, \widehat{S}_n)$. 
The principle is as follows:
\begin{itemize} 
\item At time $0$, with the values $\textbf{\em{Z}}_0 = \textbf{\em{z}}$ and $S_0 = 0$, calculate  a first 
date $R_1$ which depends on $\textbf{\em{Z}}_0$, $S_0$ and on the value that has realized 
the maximum in the calculation of $\widehat{L}_1(\widehat{v}_1,g)$.
\item Then the process is allowed to run normally until the time 
$\min\{R_1, T_1\}$. If $R_1$ comes first, it is the date of maintenance, 
if $T_1$ (the date of the first failure) comes first, the calculation is reset.
\item At time $T_1$, with the values of $\textbf{\em{Z}}_1$ and $S_1$, calculate  the second 
date $R_2$ which depends on $\textbf{\em{Z}}_1$ and $S_1$ and on the  
the value that has realized the maximum in the calculation of 
$\widehat{L}_2(\widehat{v}_2,g)$.
\item Then the process is allowed to run normally until the time 
$\min\{(T_1 + R_2),T_2\}$. If $T_1 + R_2$ comes first, it is the date of 
 maintenance, if $T_2$ comes first, reset the calculation, 
and so on until the $N$th  jump time or a total running time of $1000\ h$, whichever comes first, where maintenance will be performed 
if it has not occurred before.
\end{itemize} 
The quality of this approximation has been proved by comparing the 
expectation of the cost function of the process stopped by the above strategy 
to the true value function. This result, its proof and the precise construction 
of our stopping time procedure can be found in \cite{AAP10}.

This stopping strategy is interesting for several reasons. First, this is a 
real 
stopping time for the original PDMP which is a very strong theoretical result. Second, it requires no additional 
computation compared to those made to approximate the value function. 
This procedure can be easily performed in real time. Moreover, even if the 
original problem is an optimization {\em on average}, this stopping rule is 
path-wise and is updated when new data arrive on the history of the process at each new component failure.
%
\section{Numerical results}
\label{section:res}
The numerical procedure described above is valid for a wide class of PDMP's, see \cite{AAP10} for details. It has nice convergence properties. However, to implement it in practice, one must carefully choose the various parameters
\begin{itemize}
\item the computation horizon $N$,
\item the possibly state dependent time discretization steps $\Delta(\textbf{\em{z}})$,
\item the number of points in the quantization grids,
\item the parameters required to build quantization grids.
\end{itemize}
These parameters may not be easy to choose for a given application. One important contribution of this paper is to show that these parameters can be suitably chosen to optimize the tank.

The numerical procedure described above has been implemented on the example of the heated holdup tank. We used the exact C++ simulator of trajectories developed for \cite{ZDDG08}, suitably modified to take into account the possible failures of the command and interfaced with a matlab code for the optimization procedure. The jump horizon $N$ was empirically set to $26$ jumps, thus allowing all the trajectories to reach one of the top events  $h<4\ m$, $h>10\ m$, $\theta>100^\circ C$ or $t>1000$ hours.
%
\subsection{Quantization grids}
We encountered a new difficulty when deriving the quantization grids, due to the high cardinality of the possible modes and possibly low probability of reaching some of them.
Our mathematical model for the dynamics of the tank is hybrid: there is a discrete mode variable (the positions of the components and state of the control unit) and a continuous variable (liquid level, temperature, running time). Of course, one needs not discretize the mode variable as it can already take only finitely many values. Our procedure requires one discretization grid at each jump time of the process. However, at a given jump time, several modes can appear. For instance, at time $0$, the starting mode is $(ON, OFF, ON, 1)$. After the first jump time, one of the components has failed, so there are now 6 possible modes: $(Stuck\ ON, OFF, ON, 1)$, $(stuck\   OFF, OFF, ON, 1)$, $(ON, stuck\ ON, ON, 1)$, $(ON, stuck\ OFF, ON, 1)$, $(ON, OFF, stuck\ ON, 1)$ or $(ON, OFF, stuck\ OFF, 1)$. Table~\ref{tab2} gives the theoretical number of possible modes at each time step as well as the observed one for $3\cdot10^{9}$ simulated trajectories. After 5 jumps, the theoretical number of modes is constant equal to 18, but all the 18 modes can actually be observed only as long as the controller unit does not fail. As the probability for the controller to remain in its operational state decreases with the number of trials of the control laws, the 18 modes become increasingly rare and by the $25$-th jump time are not observed anymore, which means that the system has reached one of the top events and stopped.
\begin{table}
\begin{center}
\begin{tabular}{|c|c|c||c|c|c|}
\hline
Time&Theory&Simulations&Time&Theory&Simulation\\
\hline
n=0&1&1&n=20&18&17\\
\hline
n=1&6&6&n=21&18&16\\
\hline
n=2&18&18&n=22&18&14\\
\hline
n=3&30&30&n=23&18&7\\
\hline
n=4&25&25&n=24&18&1\\
\hline
n=5&&&n=25&18&0\\
\cline{4-6}
to n=19&18&18&n=26&18&0\\
\hline
\end{tabular}
\caption{Theoretical and observed number of modes at each time step for $3\cdot10^{9}$ Monte Carlo simulations.}\label{tab2}
\end{center}
\end{table}

The comparatively rare events are problematic for the implementation of the quantization algorithm. Indeed, one first chooses the number $k$ of discretization points then usually initializes the algorithm by throwing $k$ trajectories of the process at random. Thus, some rare modes may not be reached by the initial simulations, and the algorithm will not perform well when new trajectories are thrown reaching these modes. Indeed, the algorithm is based on a nearest neighbor search, within the nodes having the same mode as the original point. When no such mode is present, the algorithm provides highly unsatisfactory results. Therefore, we had to find a way to ensure that the initializing grids have at least one point in each observed mode for each time step. To do so, at each time step we allocated the $k$ points to the possible modes proportionally to their frequency (computed with $3\cdot10^{9}$ Monte Carlo simulations) and forcing 1 point for the modes with frequency less than $1/k$.
%
\subsection{Optimal performance and maintenance date}
We ran our optimization procedure for several number of discretization points in the quantization grids. The results we obtained are given in Table~\ref{tab3}. The optimal performance is our approximation of the value function $v$ at the starting equilibrium point whereas the last column gives the mean performance achieved by our stopping rule (obtained by $10^5$ Monte Carlo simulations).
\begin{table}
\begin{center}
\begin{tabular} {|c|c|c|c|c|c|} \hline
Number of points               & Optimal performance   &  Stopping rule            \\  \hline
200        & 334.34  &  305.55       \\  \hline
300       & 333.04  &  319.45       \\  \hline
400       & 332.95  &  322.20        \\  \hline
800       & 330.43  &  323.63        \\  \hline
1000     & 330.87  &  324.04       \\  \hline
\end{tabular} 
\caption{Optimal performance}\label{tab3}
\end{center}
\end{table} 
One can observe the convergence as the number of points in the quantization grids increases, and one can also see that our stopping rule is indeed close enough to optimality.

\begin{figure}[htb]
\begin{center}
\subfigure[distribution]
{\includegraphics[width=0.49\linewidth]{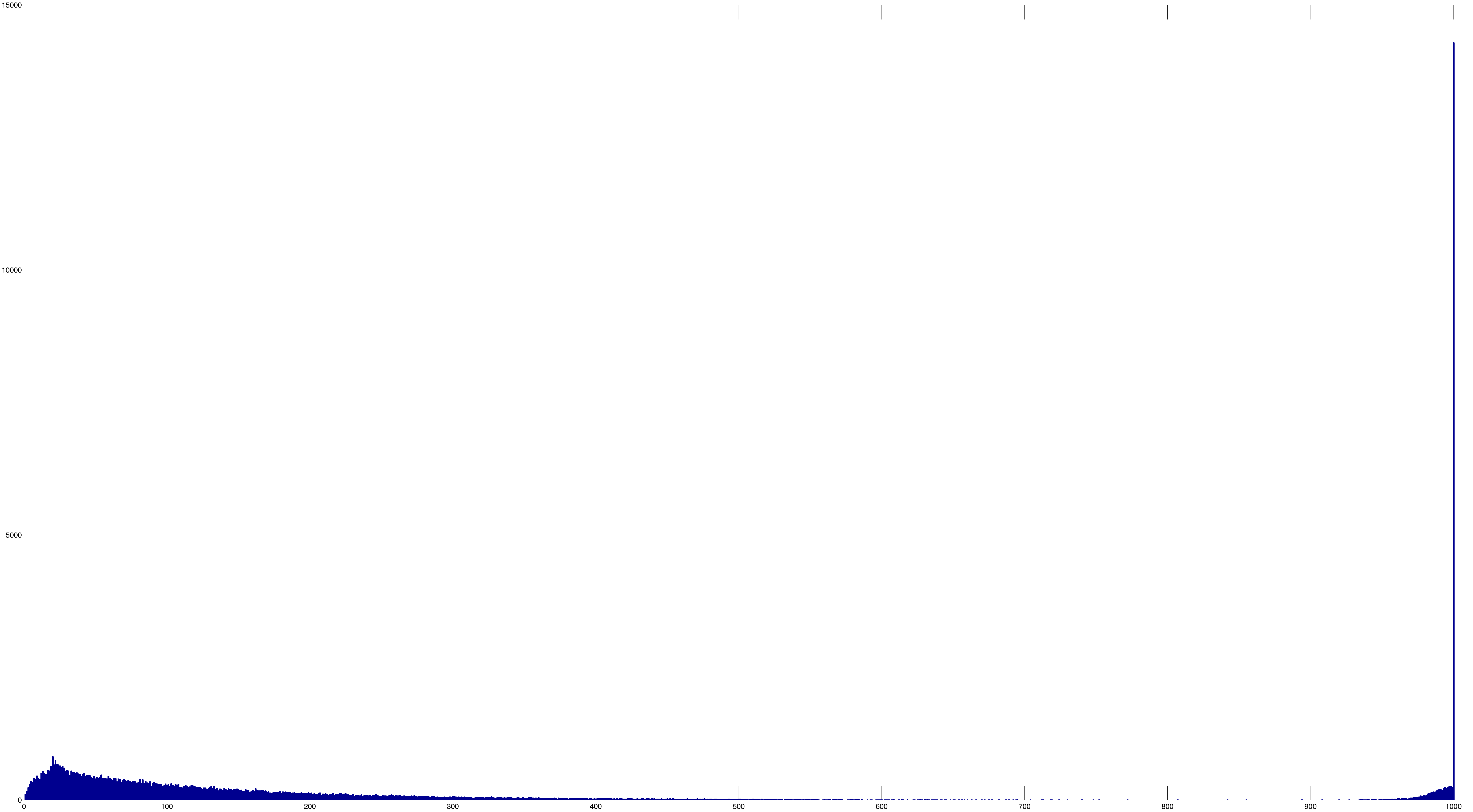}}
\subfigure[zoom of the distribution]
{\includegraphics[width=0.49\linewidth]{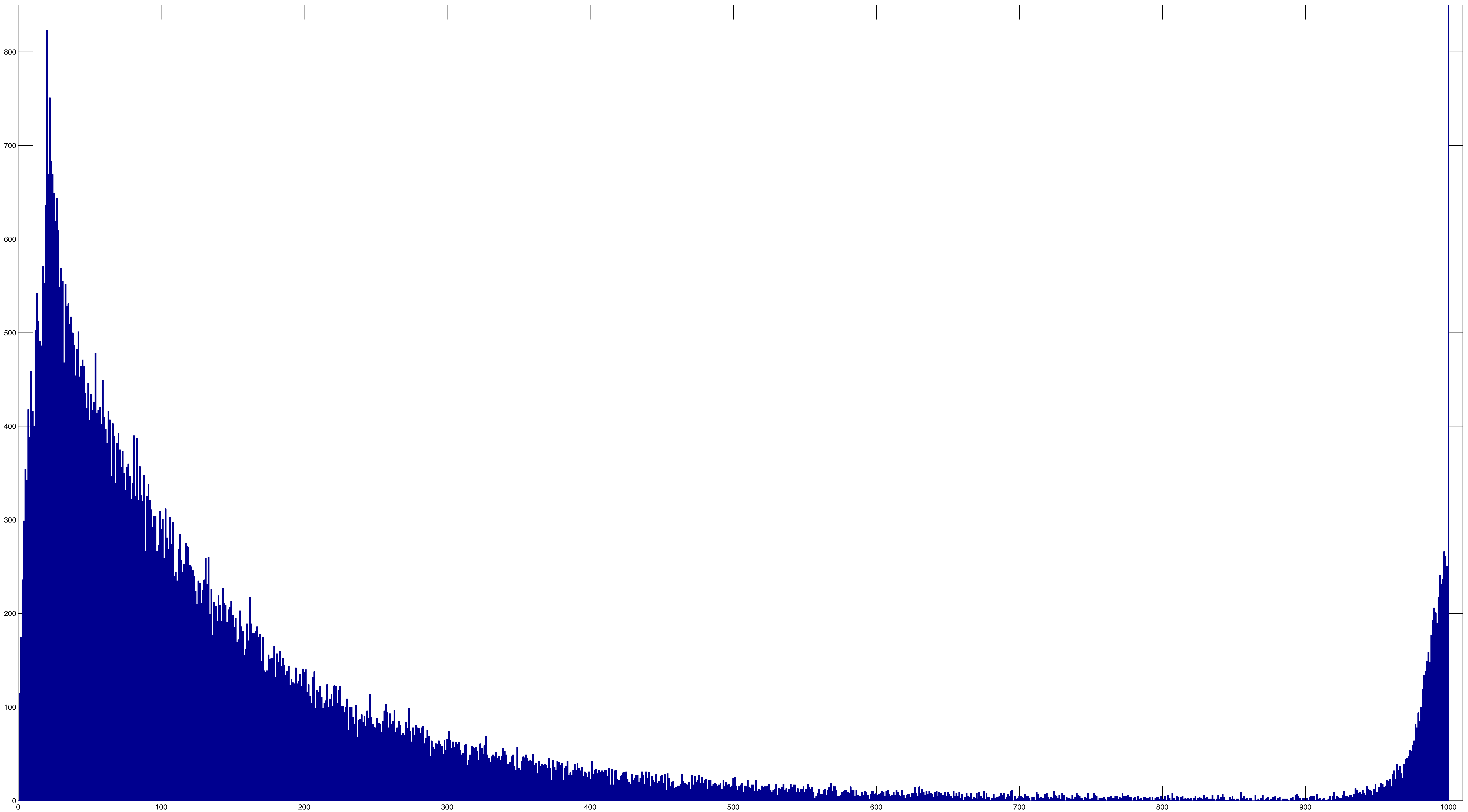}}
\caption{Histogram of the distribution of the computed maintenance time (with and without the mode at 1000)}\label{fig8}
\end{center}
\end{figure}

\begin{figure}[h]
\centering
\subfigure[Liquid level]
{\includegraphics[width=0.49\linewidth]{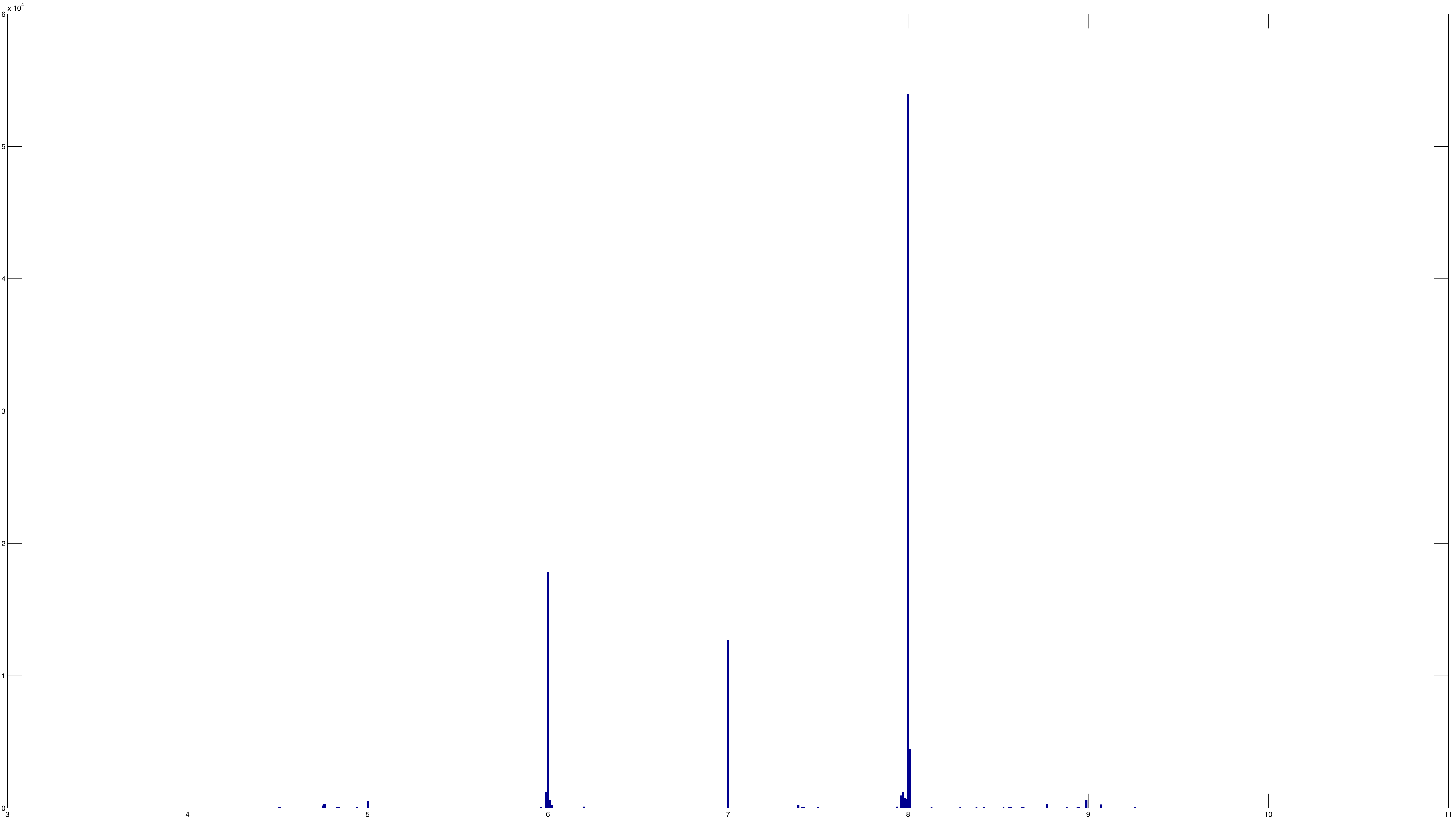}} 
\subfigure[Temperature]
{\includegraphics[width=0.49\linewidth]{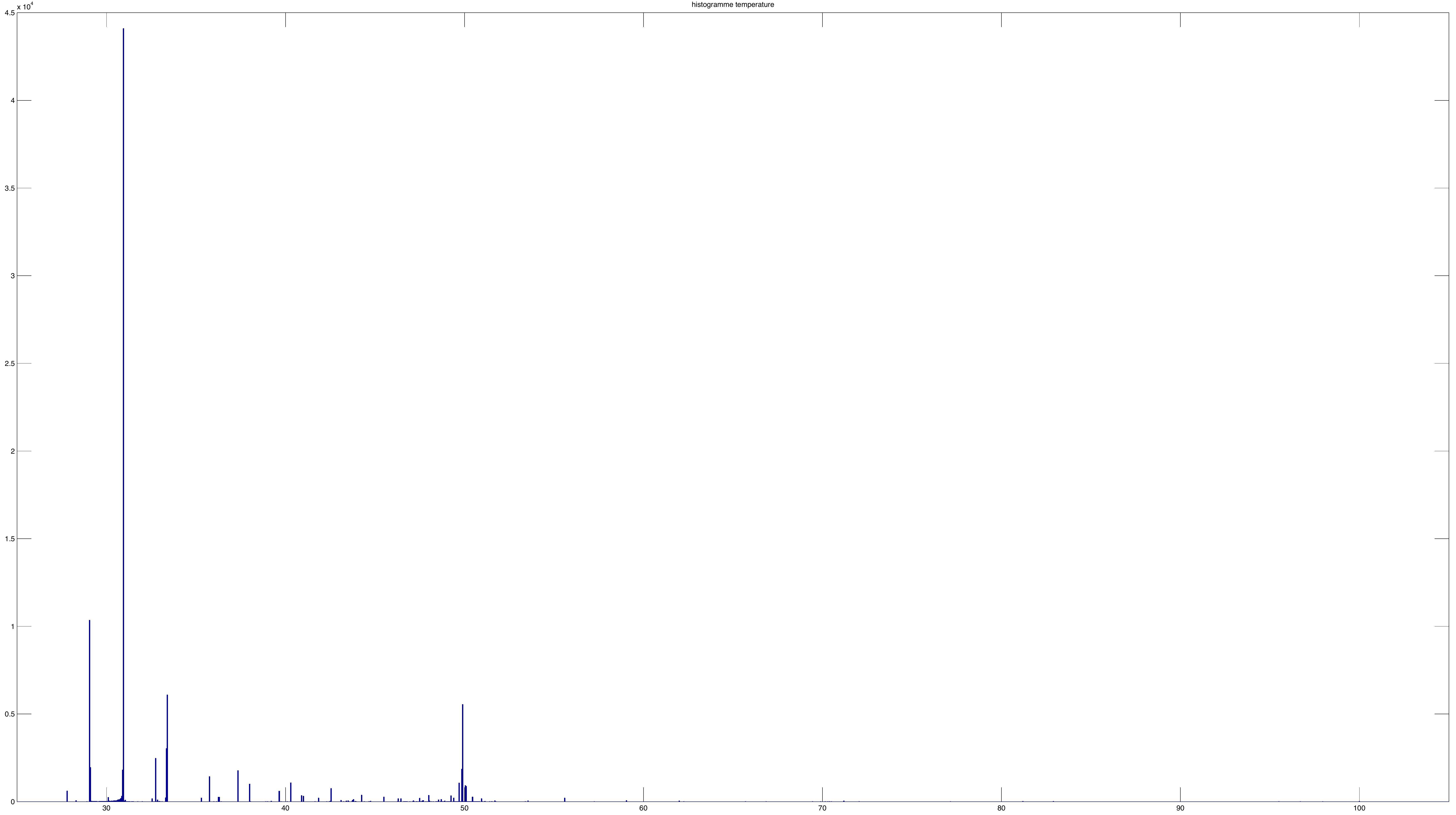}} 
\caption{Histogram of the liquid level (left) and temperature (right) at the maintenance time}
\label{fig9}
\end{figure}

We can also obtain the distributions of the maintenance time by Monte Carlo simulations ($10^5$ simulations). It is given in Figure~\ref{fig8}. The distribution is roughly bimodal, with a very high mode (14.16\% of trajectories) at time $1000$, which means that the tank remained in an acceptable state until the end of the experiment. For better readability of the histogram, only the value of $\tau^*<1000$ are plotted on the right-hand side figure.
This distribution illustrates that an average stopping rule would be far from optimal for the tank. The distribution of the liquid level and temperature at the maintenance time are given on Figure~\ref{fig9}.
The distribution of the liquid level is almost discrete with three main values at levels $6\ m$ (15.16\%), $7\ m$ (12.68\%) and $8\ m$ (45.45\%). This is natural as the liquid level is often constant, and moves very fast between these 3 states. The distribution of the temperature is also almost discrete as it is a function of the liquid level for most modes. It has a maximum at the equilibrium temperature (15,60\%). Note that the top events are never reached for the liquid level, and only $0.02\%$ of trajectories ended at $\theta=100^\circ C$. This is a strong point in favor of our procedure.
\subsection{Validation}
There is no analytical solution to our optimization problem, therefore it is impossible to compare our results with the true value function. However, it must be stressed out that there is a theoretical proof in \cite{AAP10} that this procedure converges to the true value as the number of points in the quantization grids goes to infinity.

\begin{figure}[htb]
\begin{center}
\includegraphics[width=13cm]{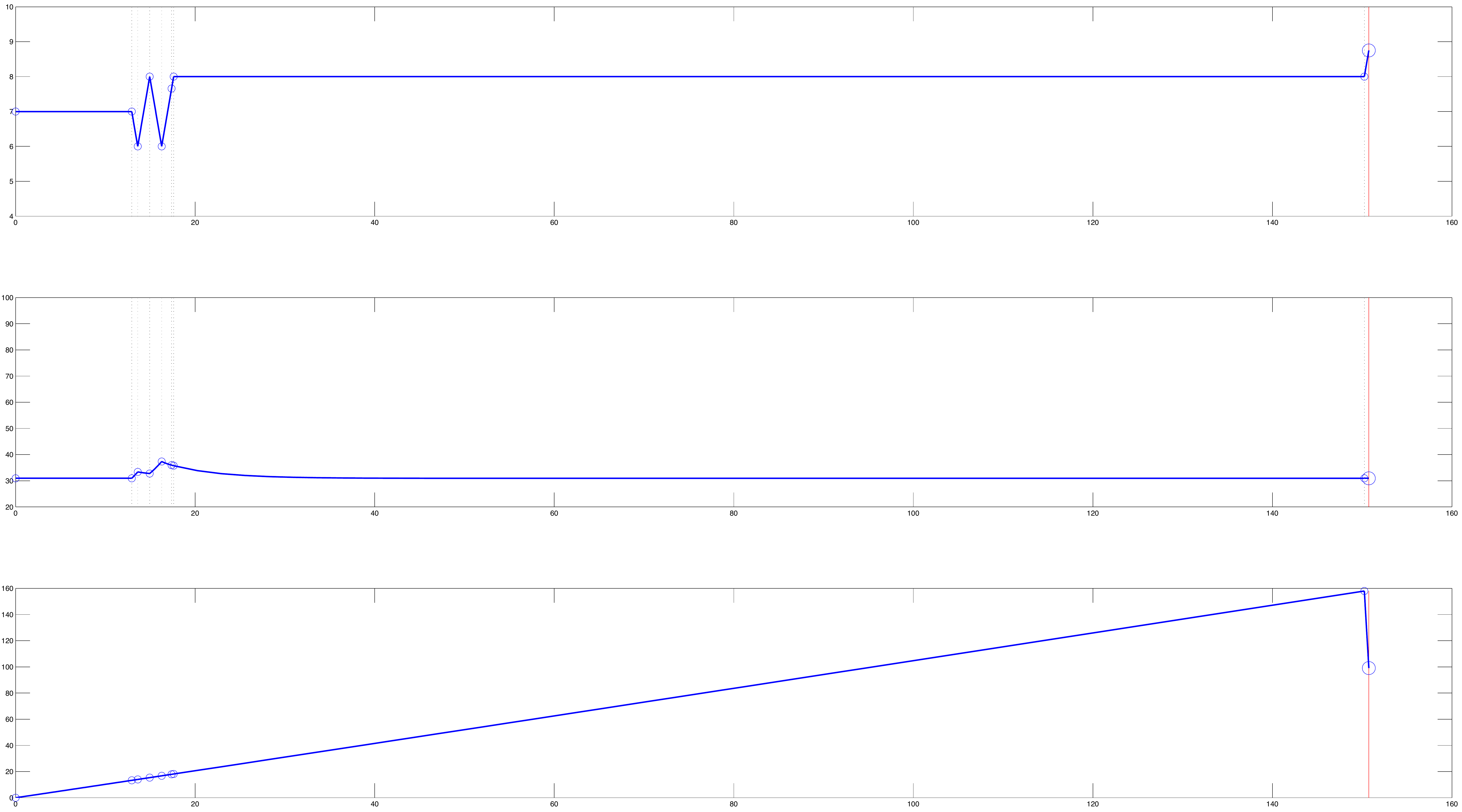}
\caption{An example of trajectory. The upper chart shows the liquid level, the middle chart the temperature, and the bottom chart the corresponding reward, as a function of time. The small circles represent the jumps of regime, and the large circle and red line the computed maintenance time.}\label{fig10}
\end{center}
\end{figure}

We can also conduct two simple kinds of experiments to validate our results. First, we can simulate new trajectories and check one by one if the intervention took place at a reasonable time. Figure \ref{fig10} shows such an interesting example. 
\begin{itemize}
\item The system starts in the equilibrium state, so that the reward function grows roughly linearly with time. 
\item At time $12.94\ h$, a jump occurs. The first unit is stuck OFF. The liquid level drops fast and the temperature rises slowly. As the liquid level and temperature are still within the acceptable bounds, the reward function is still growing roughly linearly with time. 
\item At time $13.61\ h$, the liquid level reaches the boundary $6\ m$ and the controller switches the 3 units to stuck OFF, ON and OFF. The liquid level now rises whereas the temperature drops.
\item At time $14.94\ h$, the liquid level reaches the boundary $8\ m$ and the controller switches the 3 units back to stuck OFF, OFF and ON. The liquid level drops again and the temperature rises.
\item At time $16.27\ h$, the liquid level reaches again the boundary $6\ m$ and the controller switches the 3 units back to stuck OFF, ON and OFF. The liquid level rises and the temperature drops.
\item At time $17.38\ h$, before the liquid level reaches again $8\ m$, the second unit fails and is now stuck ON. The liquid level still rises and the temperature drops.
\item At time $17.60\ h$, the liquid level reaches the boundary $8\ m$ and the controller switches the 3 units to stuck OFF, Stuck ON and ON. The liquid level now remains constant as the temperature drops down to its equilibrium state.
\item At time $150.24\ h$, the third unit is stuck OFF. The liquid level rises again rapidly and goes above the acceptable threshold. The temperature remains constant. The algorithm decides to perform a maintenance at time $151.58$. The final reward is $99.07$.
\end{itemize}

\begin{figure}[htb]
\begin{center}
\includegraphics[width=13cm]{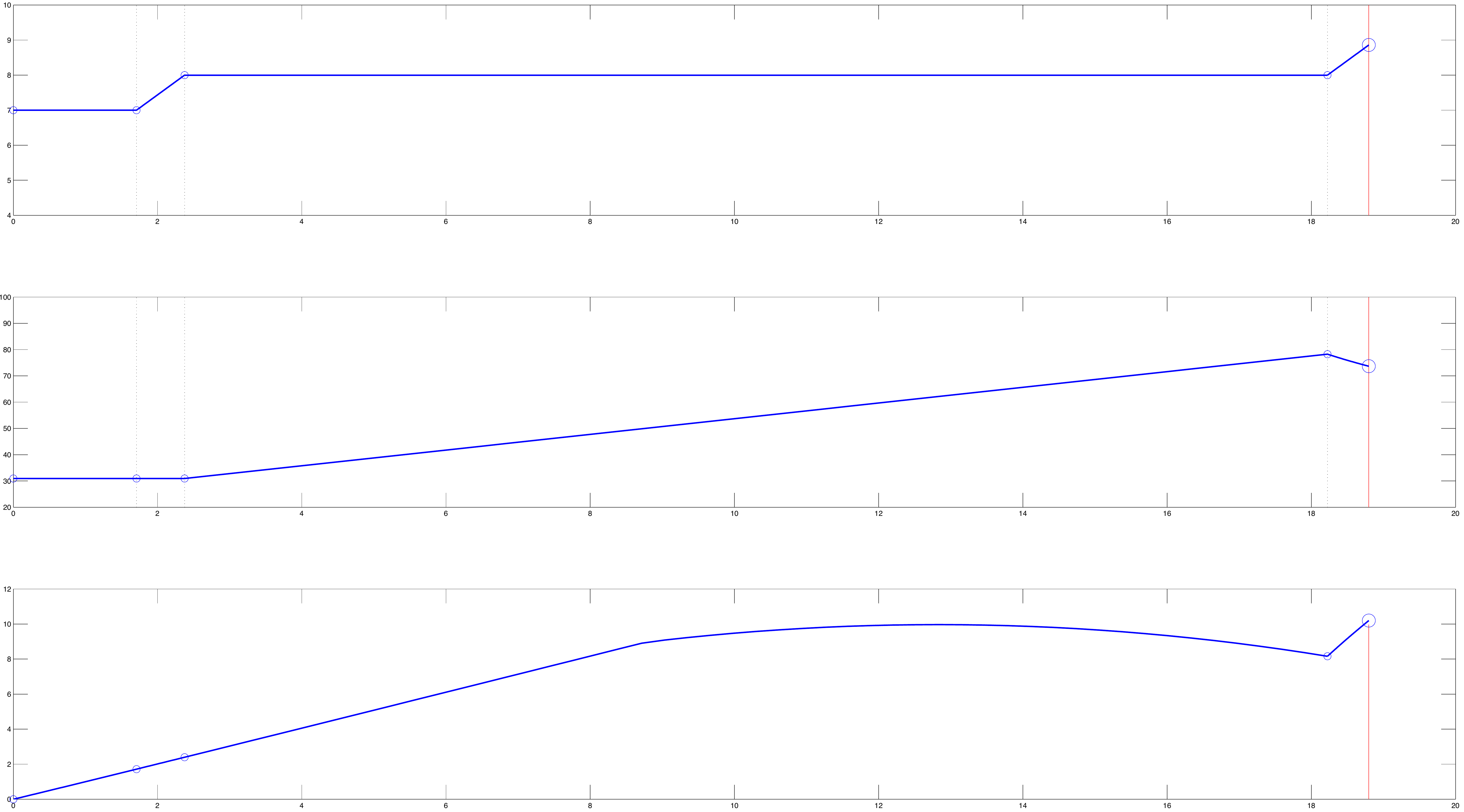}
\caption{An example of trajectory. The upper chart shows the fluid level, the middle chart the temperature, and the bottom chart the corresponding reward, as a function of time. The small circles represent the jumps of regime, and the large circle and red line the computed maintenance time.}\label{fig11}
\end{center}
\end{figure}

Another example is given in Figure~\ref{fig11}.
\begin{itemize}
\item The system starts in the equilibrium state, so that the reward function grows linearly with time. 
\item At time $1.71\ h$, a jump occurs. The third unit is stuck OFF. The liquid level rises and the temperature remains constant. As the liquid level and temperature are still within the acceptable bounds, the reward function is still growing linearly with time. 
\item At time $2.37\ h$, the liquid level reaches the threshold 8m and triggers the command. The solicitation of the command is successful, and the first unit is turned OFF. The current state is now OFF for units 1 and 2, and stuck OFF for unit 3. The liquid level remains constant at $8\ m$, but the temperature rises. At about $9\ h$, the temperature crosses the threshold $50^\circ C$ so that the reward function is now slowly decreasing. Note that this does not trigger an immediate maintenance. 
\item At time $18.22\ h$, another jump occurs and the second unit is now stuck ON, causing the liquid level to rise again, but the temperature to decrease. As a result, the reward function is now increasing again. 
\item The algorithm then selects the maintenance time to be $18.8\ h$, before the liquid level reaches the level $10\ m$ causing a total failure of the system. 
\end{itemize}
Note that at this intervention time, the first unit is not stuck, so that other jumps may happen in the future. In this special example, the intervention time occurs when the reward is maximal.

Second, one can simply compare the performances of the system when no maintenance is performed to those with our maintenance rule. The results are summarized in Table~\ref{tab4}. These results are obtained with $10^5$ Monte Carlo simulations and illustrate the power of our procedure as the mean performance is increased by $156\%$ and the top events are almost always avoided.
\begin{table}
\begin{center}
\begin{tabular}{|l|c|c|}
\hline
&without maintenance&with maintenance\\
\hline
mean performance&211.80&{330.87}\\
\hline
null gain&80.33\%&{0.02\%}\\
\hline
$6\leq h\leq 8$&28.25\%&{90.02\%}\\
\hline
$\theta\leq 50^\circ C$&80.33\%&{95.09\%}\\
\hline
$h=4$&16.65\%&0\%\\
\hline
$h=10$&54.55\%&0\%\\
\hline
$\theta=100$&9.13\%&0.02\%\\
\hline
\end{tabular}
\caption{Comparison of the performances with and without maintenance}\label{tab4}
\end{center}
\end{table}
%
\section{Conclusion}
\label{section:conclusion}
The numerical method described in \cite{AAP10} has been applied to a well known test case of dynamic reliability to approximate the value function of the optimal stopping problem and an $\epsilon$-optimal stopping time for a piecewise-deterministic Markov process, that is the maintenance date for the tank. The quantization method proposed can sometimes be costly in computing time, but has a very interesting property: it can be calculated off-line. Moreover it depends only on the dynamics of the model, and not on the cost function chosen, or the actual trajectory of the specific process one wants to monitor. The calculation of the optimal maintenance time is done in real time, making our procedure applicable in practice. The optimal maintenance time is updated at the changes of mode and has a conditional threshold form, which allows scheduling maintenance services in advance.

If one only changes the reward function $g$ without changing the dynamics of the tank, one just has to run the optimization part of the algorithm, and not the quantization grids. This can be done in real time. If one wants to change the dynamics of the system, or add some components, one has to rewrite the simulation code for the system, and with this new code re-run the quantization grids, which can be quite long. However, the general methodology is valid for a wide class of piecewise deterministic Markov processes and not at all specific to the tank.

The method has been implemented on the heated hold-up tank. The main characteristic of this system is that it can be modeled by a stochastic hybrid process, where the discrete and continuous parts interact in a closed loop. The optimization problem under study has no analytic solution. However, our method is based on a rigorous mathematical construction with proof of convergence. In addition, simple comparisons between no motoring and our policy also prove its practical validity with a significant improvement of the performance of the system (the mean performance is increased by $156\%$ and the top events are almost always avoided).

Our next project is to extend this research in two main directions. First, we could allow only partial repair of the system. The problem will then be to find simultaneously the optimal times of maintenance and optimal repair levels. Mathematically, it is an impulse control problem, which complexity exceeds widely that of the optimal stopping. Second, our method requires a perfect observation of the state process at the jump times. It would be interesting to extend our results to a noisy observation of the process, as often happens in real life. 

\section*{Acknowledgements}
This work was partially funded by the French National Research, ANR grant FAUTOCOES, program ARPEGE, ANR-09-SEGI-004. The authors also thank Pr. Fran\c cois Dufour for initial discussions on how to define an interesting optimal stopping problem for the test-case of the tank. This paper was presented at ESREL 2012 conference.


\begin{thebibliography}{10}

\bibitem{aldemir87a}
{\sc Aldemir, T.}
\newblock Computer-assisted markov failure modeling of process control systems.
\newblock {\em IEEE Transactions on Reliability 36}, 4 (1987), 133--144.

\bibitem{AJ99}
{\sc Aven, T., and Jensen, U.}
\newblock {\em Stochastic models in reliability}, vol.~41 of {\em Applications
  of Mathematics}.
\newblock Springer, New York, 1999.

\bibitem{cojazzi1996a}
{\sc Cojazzi, G.}
\newblock The dylam approach for dynamic reliability analysis of systems.
\newblock {\em Reliability Engineering and System Safety 52\/} (1996),
  279--296.

\bibitem{Davis93}
{\sc Davis, M. H.~A.}
\newblock {\em Markov models and optimization}, vol.~49 of {\em Monographs on
  Statistics and Applied Probability}.
\newblock Chapman \& Hall, London, 1993.

\bibitem{AAP10}
{\sc de~Saporta, B., Dufour, F., and Gonzalez, K.}
\newblock Numerical method for optimal stopping of piecewise deterministic
  {M}arkov processes.
\newblock {\em Ann. Appl. Probab. 20}, 5 (2010), 1607--1637.

\bibitem{JRR12}
{\sc de~{S}aporta, B., {D}ufour, F., {Z}hang, H., and {E}legbede, C.}
\newblock {O}ptimal stopping for the predictive maintenance of a structure
  subject to corrosion.
\newblock {\em Proc. of the Inst. of Mechanical Engineers, Part O: Journal of
  Risk and Reliability 226}, 2 (2012), 169--181.

\bibitem{devooght97}
{\sc Devooght, J.}
\newblock {\em Dynamic reliability}.
\newblock Advances in nuclear science and technology. Chapman and Hall, Berlin,
  1997.

\bibitem{DD02}
{\sc Dufour, F., and Dutuit, Y.}
\newblock Dynamic reliability: A new model.
\newblock In {\em Proceedings of ESREL 2002 Lambda-Mu 13 Conference\/} (2002),
  pp.~350--353.

\bibitem{dutuit1997a}
{\sc Dutuit, Y., Ch\^atelet, E., Signoret, J., and Thomas, P.}
\newblock Dependability modelling and evaluation by using stochastic petri
  nets: application to two test cases.
\newblock {\em Reliability Engineering and System Safety 55\/} (1997),
  117--124.

\bibitem{grall}
{\sc Grall, A., B\'erenguer, C., and Dieulle, L.}
\newblock A conditional-based maintenance policy for a stochastically
  deteriorating system.
\newblock {\em Reliability Engineering and System Safety 76}, 2 (2002),
  167--180.

\bibitem{gugerli86}
{\sc Gugerli, U.~S.}
\newblock Optimal stopping of a piecewise-deterministic {M}arkov process.
\newblock {\em Stochastics 19}, 4 (1986), 221--236.

\bibitem{li2011a}
{\sc Li, J., Moslehb, A., and Kanga, R.}
\newblock Likelihood ratio gradient estimation for dynamic reliability
  applications.
\newblock {\em Reliability Engineering and System Safety 96\/} (2011),
  1667--1679.

\bibitem{marseguerra1994a}
{\sc Marseguerra, M.}
\newblock Approximated physical modelling in dynamic psa using artificial
  neural networks.
\newblock {\em Reliability Engineering \& System Safety 45\/} (1994), 47--56.

\bibitem{MZ95}
{\sc Marseguerra, M., and Zio, E.}
\newblock The cell-to-boundary method in {M}onte {C}arlo-based dynamic {PSA}.
\newblock {\em Reliability Engineering and System Safety 48}, 3 (1995),
  199--204.

\bibitem{MZ96}
{\sc Marseguerra, M., and Zio, E.}
\newblock {M}onte {C}arlo approach to {PSA} for dynamic process systems.
\newblock {\em Reliability Engineering and System Safety 52}, 3 SPEC. ISS.
  (1996), 227--241.

\bibitem{pages98}
{\sc Pag{\`e}s, G.}
\newblock A space quantization method for numerical integration.
\newblock {\em J. Comput. Appl. Math. 89}, 1 (1998), 1--38.

\bibitem{pages05}
{\sc Pag{\`e}s, G., and Pham, H.}
\newblock Optimal quantization methods for nonlinear filtering with
  discrete-time observations.
\newblock {\em Bernoulli 11}, 5 (2005), 893--932.

\bibitem{pages04b}
{\sc Pag{\`e}s, G., Pham, H., and Printems, J.}
\newblock An optimal {M}arkovian quantization algorithm for multi-dimensional
  stochastic control problems.
\newblock {\em Stoch. Dyn. 4}, 4 (2004), 501--545.

\bibitem{pages04a}
{\sc Pag{\`e}s, G., Pham, H., and Printems, J.}
\newblock Optimal quantization methods and applications to numerical problems
  in finance.
\newblock In {\em Handbook of computational and numerical methods in finance}.
  Birkh\"auser Boston, Boston, MA, 2004, pp.~253--297.

\bibitem{schoenig2006a}
{\sc Schoenig, R., Aubry, J., Cambois, T., and Hutinet, T.}
\newblock An aggregation method of markov graphs for the reliability analysis
  of hybrid systems.
\newblock {\em Reliability Engineering and System Safety 91\/} (2006),
  137--148.

\bibitem{siu1994a}
{\sc Siu, N.}
\newblock Risk assessment for dynamic systems: an overview.
\newblock {\em Reliability Engineering and System Safety 43\/} (1994), 43--73.

\bibitem{tombuyses1996a}
{\sc Tombuyses, B., and Aldemir, T.}
\newblock Continous cell-to-cell mapping and dynamic psa.
\newblock In {\em Proceedings of ICONE 4 conference\/} (1996), pp.~431--438.

\bibitem{vannoortwijk}
{\sc van Noortwijk~J.M.}
\newblock A survey of the application of gamma processes in maintenance.
\newblock {\em Reliability Engineering and System Safety 94}, 1 (2009), 2--21.

\bibitem{ZDDG08}
{\sc {Z}hang, H., {D}ufour, F., Dutuit, Y., and Gonzalez, K.}
\newblock Piecewise deterministic markov processes and dynamic reliability.
\newblock {\em Proc. of the Inst. of Mechanical Engineers, Part O: Journal of
  Risk and Reliability 222}, 4 (2008), 545--551.

\end{thebibliography}
\end{document}